\documentclass{article}
\usepackage{enumerate}
\usepackage{amsmath}
\usepackage{graphicx}
\usepackage{subfig}
\usepackage{tensor}
\usepackage{float}
\begin{document}
\title{A Mathematical Random Number Generator (MRNG)}
 \author{
Osvaldo Skliar\thanks{Escuela de Inform\'atica, Universidad Nacional, Costa Rica.  E-mail: oskliar@costarri\-cense.cr}
\and Ricardo E. Monge\thanks{Escuela de Inform\'atica, Universidad Nacional, Costa Rica. E-mail: ricardo.monge.gap\-per@una.cr}
\and Sherry Gapper\thanks{Universidad Nacional, Costa Rica. E-mail: sherry.gapper.morrow@una.cr}
\and Guillermo Oviedo\thanks{Universidad Latina, San Pedro, Costa Rica. E-mail: oviedogmo@gmail.com}}

\maketitle
\begin{abstract}
A  novel Mathematical Random Number Generator (MRNG) is presented here. In this case, "mathematical" refers to the fact that to construct that generator it is not necessary to resort to a physical phenomenon, such as the thermal noise of an electronic device, but rather to a mathematical procedure. The MRNG generates binary strings -- in principle, as long as desired -- which may be considered genuinely random in the sense that they pass the statistical tests currently accepted to evaluate the randomness of those strings. From those strings, the MRNG also generates random numbers expressed in base 10.  An MRNG has been installed as a facility on the following web page: www.appliedmathgroup.org. This generator may be used for applications in tasks in: a) computational simulation of probabilistic-type systems, and b) the random selection of samples of different populations. Users interested in applications in cryptography can build another MRNG, but they would have to withhold information -- specified in section 5 -- from people who are not authorized to decode messages encrypted using that resource.

\end{abstract}



\section{Introduction}

There are diverse characterizations of the concept of randomness. The notions of random numbers and pseudo-random numbers are closely related to the concept of randomness. Therefore, it is necessary to specify which characterization of the concept of randomness is applied when referring to this type of numbers. Notions of randomness different from that used in this paper will be discussed in sections 4 and 5. 

The present section is devoted mainly to presenting the objective of this article. To this end, a preliminary specification will be provided of the notion of randomness used here; that is, the one usually accepted in literature on statistics when ``random numbers'' are mentioned. That characterization will be completed in sections 2 and 3. 

In \cite{bhrng} the notion of an ideal random-bit generator (IRBG) was introduced in the following way:

``For a clearer understanding of the methods presented here, it is useful to consider an ideal machine that every so often (for instance, every millisecond) generates a digit with the following two characteristics: a) the digit is either a zero (0) or a one (1), and b) there is the same probability that the digit will be a 0 or a 1 (that is, for each digit generated by this ideal machine, the probability that it could be a 1 is 0.5, and that it could be a 0 is 0.5). This ideal machine will be called an `ideal random-bit generator' (IRBG). It will also be supposed that it has an ideal `memory' device (e.g., an `ideal hard disk') which makes it possible to store a sequence of bits (generated by the IRBG) that is as long, or that contains as many bits, as desired. This sequence will be called a `binary string'.'' (p. 267)

To come as near as possible to the IRBG and for the highest quality, it would appear to be indispensable to resort to a physical process, such as atmospheric noise, thermal noise in electronic systems or radioactive decay, in which probabilistic-type events are accepted to be present inevitably. Nevertheless, the use of a physical device may be unnecessary if attention is not given to the IRBG itself, but rather to its product, its output: a binary string with the previously specified characteristics. Suppose that the scientific community is presented a deterministic-type procedure which makes it possible to obtain a binary string made up of as many bits as desired. Suppose also that for someone who is not too familiar with the details of the procedure used to generate them, each one of these strings cannot be distinguished from the binary string generated by an IRBG. How can it be determined whether these binary strings actually do have the property which they are purported to have? The answer is evident: Statistical tests can be applied just as they are commonly used to inquire about the quality of numerical sequences produced by a) physical generators of supposedly random numbers, and by b) generators of pseudo-random numbers, usually consisting of programs based on certain algorithms.

The deterministic procedure presented here is of an exclusively mathematical character; that is, it does not resort to any physical type of resource at all. It generates a binary string which is as long as desired and has the following characteristic: if that binary string is subjected to all of the accepted statistical tests used to determine the randomness of a binary string, it passes those tests. In other words, according to the accepted statistical criteria used to establish that a binary string is random, that string must be considered random.

Naturally, whoever is familiar in detail (what is meant by the expression ``in detail'' will be indicated with no ambiguity whatsoever in Section 5) with the procedure used to generate this binary string can compute, for every digit of the string, its numeric value (that is, determine whether each digit is either a 0 (zero) or a 1 (one)).  This occurs because the string is generated by an entirely deterministic algorithmic procedure. However, for someone who is not familiar in detail with that procedure,  the numeric value of any digit of the binary string is absolutely unforeseeable, regardless of the amount of preceding digits whose value has already been known.  Suppose, for example that this person is given the numeric values of the $10^{12}$ digits preceding a given digit. In that case, no matter how vast this person's mathematical knowledge is and no matter what computational tools he has, he will only be able to state that there is a probability of $\frac{1}{2}$ that the digit is a 0 and the same probability of  $\frac{1}{2}$ that the digit is a 1.

Once the binary string is available, it is simple, as will be explained below, to obtain from it genuine random numbers (not just pseudo-random numbers), made up of a sequence of digits such that each of them may be one of the following: 0, 1, 2, 3, 4, 5, 6, 7, 8, 9 (i.e., the digits corresponding to numerical base 10). 

The objective of this paper is to present this mathematically-based generator of genuine random numbers obtained from the binary string introduced above.

\section {A mathematical random number generator \\(MRNG)}

Roots (whose degrees are prime numbers) of prime numbers were used to obtain the random binary string mentioned in section 1. Consider, for example, the following roots of the prime numbers 5 and 17:
\begin{equation*}
5^{\frac{1}{3}} \quad  \textrm{and}\quad 17^{\frac{1}{3}}
\end{equation*}

The computational techniques currently available make it possible to obtain these cubic roots (which are irrational numbers) having, for example, hundreds of thousands of decimal digits. GMP \cite{bgmp} and MPFR \cite{bmpfr} were used to compute mathematically correct roots.

Suppose that these roots are actually obtained with one hundred thousand decimal digits each. The whole parts of both roots are discarded and only their decimal parts will be used. Then, a certain amount, between 40 and 100, of the initial digits of each of those decimal parts is also discarded. This way of proceeding will be justified below. In the case at hand, the first 50 digits of each of the decimal parts was discarded. The corresponding sequences of digits from digit 51 to digit 70, for each of the roots, are shown in figure 1.  In addition, the one hundred thousandth digit is given for each root.

\begin{figure}[H]
\centering
\includegraphics[width=4in]{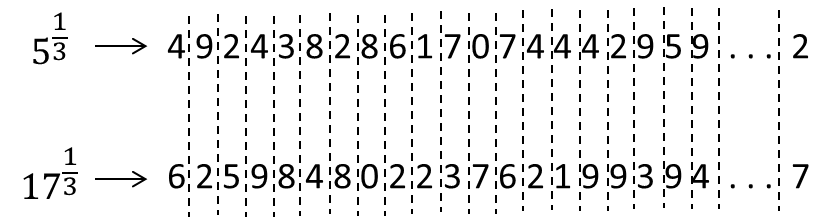}
\caption{The digits of the decimal part from digit 51 to digit 70 are shown for each of the roots $5^{\frac{1}{3}}$ and $17^{\frac{1}{3}}$.  Digit 100,000 is also given for each of these decimal parts.}
\label{f1}
\end{figure}

The first digit of those remaining (digit 51) of the decimal part of $5^{\frac{1}{3}}$ is compared to the corresponding digit 51 of the decimal part of $17^{\frac{1}{3}}$. If the first is greater than the second, it is considered to have gotten a 1 (one). If the first is less than the second, it is considered to have gotten a 0 (zero). If both digits are equal, it is considered that no digit was obtained as a result of the comparison described.

The same procedure is followed with digits $52, 53,\dots,100,000$. In this case as a result of each of the comparisons specified, the elements of the dyads appearing explicitly in figure 1 -- $(4,6), (9,2), (2,5),\dots,(9,4)$ and $(2,7)$ --  either a 0 or a 1 is generated. This, of course, is because neither of these dyads is made up of two equal elements. As explained above, dyads composed of two equal elements -- $(0,0), (1,1), \dots,(9,9)$ -- do not contribute digits to the binary string under construction.

The digits obtained by the procedure described, from the dyads appearing explicitly in figure 1, are shown in figure 2.
\begin{figure}[H]
\centering
\includegraphics[width=3in]{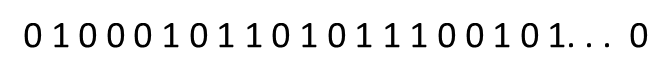}
\caption{The sequence of digits obtained by the making the comparisons described above is given here.  Only the digits obtained from the dyads appearing explicitly in figure 1 are shown.}
\label{f2}
\end{figure}

It is possible to consider that the binary string mentioned in figure 3 is the result of applying a particular operator, which will be called $O$, to the operand composed of the ordered pair $(5^{\frac{1}{3}},17^{\frac{1}{3}})$.
\begin{equation*}
O\left(5^{\frac{1}{3}},17^{\frac{1}{3}}\right)=\texttt{01}\ldots \texttt{0}\ldots \texttt{1}
\end{equation*}

In this paper, the following two sets of prime numbers were used: a) the set of the first 10,000 prime numbers, and b) the set of the next 10,000 prime numbers. In the first of these sets -- $C_{1}$ -- the prime numbers are considered in increasing order as the appear: $2, 3, 5,\dots$. In contrast, in the second of the sets -- $C_{2}$ -- the prime numbers are considered in the order that they are in after  having subjected that set of 10,000 prime numbers to any permutation of the $(10,000!-1)$ permutations such that each of them is different from the identical permutation. (``Identical permutation'' is used to refer to that in which the elements of $C_{2}$ are left in their natural increasing order.) The permutation used will be called $P_{1}$. The subscript 1 indicates that it is the first permutation in a set of permutations of $C_{2 }$, which will be carried out.

The following notation is introduced: The prime numbers  $2, 3, 5, 7,\dots$ will be called $p_{1,1},p_{2,1},p_{3,1},p_{4,1},\ldots$, respectively. The letter ``$p$'' is used to recall that each of these numbers is a prime number. The first subscript of $p$ refers to the order corresponding to each of the prime numbers in set $C_{1}$. The second subscript of $p$ -- 1 -- is used to indicate that each of those prime numbers belongs to $C_{1} $. Thus $C_{1}$ is an ordered set that can be symbolized in the following way:
\begin{equation*}
C_{1}=\left(p_{1,1},p_{2,1},\ldots,p_{10000,1}\right)
\end{equation*}

To refer to the prime numbers that are elements belonging to the ordered  set $C_{2}$, three subscripts will be used: The first of these subscripts will indicate, for each prime number belonging to $C_{2}$, the order in which it has ended up after permutation $P_{1}$; the second subscript -- 2 -- will be used to indicate that each of these prime numbers belongs to $C_{2}$; and the third subscript -- 1, in this case -- indicates that the prime numbers belonging to $C_{2}$ have ended up in the order resulting from permutation $P_{1}$. Thus the ordered set $C_{2}$ can be expressed in the following way: 
\begin{equation*}
C_{2}=\left(p_{1,2,1},p_{2,2,1},p_{3,2,1},\ldots,p_{10000,2,1}\right)
\end{equation*}

Recall that  $O((p_{1,1})^{\frac{1}{2}},(p_{1,2,1})^{\frac{1}{2}})$ is equal to a particular binary string. This string will be the first substring of the ``long binary string'' desired. The second substring, which will be placed immediately after the first, will be $O((p_{2,1})^{\frac{1}{2}},(p_{2,2,1})^{\frac{1}{2}})$. This second binary substring will be placed immediately after the first binary substring. This is how the ``long binary string'' desired will begin to be built. The third binary substring, which is placed after the binary string obtained as a result of having placed the second binary substring after the first one, is as follows: $O((p_{3,1})^{\frac{1}{2}},(p_{3,2,1})^{\frac{1}{2}})$. The next substring is $O((p_{4,1})^{\frac{1}{2}},(p_{4,2,1})^{\frac{1}2}$. This substring is placed after the existing binary string which was obtained by the procedure described above. The term ``concatenation operation'' will be used to refer to the process of placing new substrings after the existing string. In this way, substring $O((p_{10000,1}), (p_{10000,2,1}))$ is finally concatenated to the existing string. 

Once the pair of prime numbers made up of the last prime number of $C_{1}$ -- $p_{10000,1}$ -- and the last prime number of $C_{2}$ -- $p_{10000,2,1}$ -- have been used, one must proceed by leaving $C_{1}$ unmodified. On the other hand, the prime numbers belonging to $C_{2}$, are subjected to a second permutation $P_{2}$ which has the following characteristic: No prime number in $C_{2}$ will have, as a consequence of having done $P_{2}$, the same order number that it had previously. A simple way of obtaining a permutation of this type is this: The first prime number in $C_{2}$ will occupy the second place according to $P_{2}$; the second prime number in $C_{2}$ will occupy the third place according to $P_{2}$, and so on successively until the prime number corresponding order number 10,000 occupies the first place, according to $P_{2}$. In addition, the operation will continue to be carried out with the cubic roots (and not as before with square roots) of each new pair of prime numbers resulting from having carried out $P_{2}$.  Hence the new substrings are:
\begin{gather*}
O((p_{1,1})^{\frac{1}{3}},(p_{1,2,2})^{\frac{1}{3}}), O((p_{2,1})^{\frac{1}{3}},(p_{2,2,2})^{\frac{1}{3}}), O((p_{3,1})^{\frac{1}{3}},(p_{3,2,2})^{\frac{1}{3}}), \dots, \\ O((p_{10000,1})^{\frac{1}{3}},(p_{10000,2,2})^{\frac{1}{3}})
\end{gather*}

The elements in $C_{2}$ are then subjected to permutation $P_{3}$ such that no element of $C_{2}$ will have the same order number that it had as a consequence of having carried out the preceding permutations $P_{1}$  and  $P_{2}$. (In this paper to obtain $P_{3}$, use was made of the same approach of displacement of the elements in $C_{2}$ as that applied when doing  $P_{2}$.) Thus the following subsequences are obtained and duly concatenated:
\begin{gather*}
O((p_{1,1})^{\frac{1}{5}},(p_{1,2,3})^{\frac{1}{5}}), O((p_{2,1})^{\frac{1}{5}}, (p_{2,2,3})^{\frac{1}{5}}), O((p_{3,1})^{\frac{1}{5}},(p_{3,2,3})^{\frac{1}{5}}), \dots, \\ O((p_{10000,1})^{\frac{1}{5}}, (p_{10000,2,3})^{\frac{1}{5}})
\end{gather*}

The elements of $C_{2}$ are then subjected to another permutation $P_{4}$ with the same characteristic as in the previous permutations with regard to the order of the elements in $C_{2}$; and the operation is carried out using the seventh roots of the pairs belonging to prime numbers. Given the way specified in which the successive permutations $P_{2}, P_{3}, P_{4},\dots$, are obtained, no more than 10,000 permutations should be carried out so that the pairs of prime numbers used will be different from those used previously.

In order to express clearly and concisely the procedure described above, the subsequence concatenation operator $\mbox{\Large S}$ will be introduced; each of the subsequences is a specific binary string.

Consider, for example, the following binary strings $B_{1}$ and  $B_{2}$:
\begin{eqnarray*}
B_1&:& 011101 \\
B_2&:& 11101010
\end{eqnarray*}

Therefore, 
\begin{equation*}
S(B_1,B_2)=01110111101010
\end{equation*}

Just as in a summation  it is possible to specify from which term and to which other term they are added (as in $\sum\limits_{i=1}^{n} t_{i}$ ), and just as in a product of a sequence it is possible to specify from which factor and to which other factor they are multiplied (as in $\prod\limits_{i=1}^{n} f_{i}$), when using the operator $\mbox{\large S}$ for the concatenation of subsequences, it is also possible to indicate from which binary string and to which other binary string they will be concatenated: $ \overset{n}{\underset{i=1} {\mbox{\Large S}}} B_{i}$. And just as summations and products of sequences may be given a double or triple use (or even more), likewise it can be done with the $\mbox{\Large S}$ concatenation operator.
 
Consider the ordered set of prime numbers in increasing order: 

\begin{equation*}
C_{p}=(2, 3, 5, 7, 11, 13,\dots)
\end{equation*} 

For any of these prime numbers, the generic name $p_{j}$ will be used. If $j=1$, $p_{j}=2$; if $j=2$, $p_{j}=3$; If $j=3$, $p_{j}=5$, etc. Then the binary string $B$ (obtained as a consequence of having carried out 10,000 permutations of the type described, and with the operations specified above having been done with the respective pairs of prime numbers) can be characterized as follows:
\begin{equation*}
B=\overset{10,000}{\underset{j=1} {\mbox{\Large S}}} \quad \overset{10,000}{\underset{i=1} {\mbox{\Large S}}} \left(O\left((p_{i,1})^{\frac{1}{p_j}},(p_{i,2,j})^{\frac{1}{p_j}}\right)\right)
\end{equation*}

Note that $j$ assumes values from 1 to 10,000. If $j=1$, $p_{j}=2$ and $p_{1,2,j}=p_{1,2,1}$. The third subscript of $p_{1,2,1}$ indicates that the prime numbers belonging to $C_{2}$ were subjected to permutation $P_{1}$. If $j=2$, $p_{j}=3$ and $p_{1,2,j}=p_{1,2,2}$. The third subscript of $p_{1,2,2}$ indicates that the prime numbers belonging to $C_{2}$ were subjected to the permutation $P_{2}$; and so on, successively until $j$ reached the value of 10,000. For each value of $j$, 10,000 binary substrings of those considered are obtained, since $i$ varies from $i=1$ to $i=10,000$.

If one wants to lengthen the binary string obtained with the procedure described, the next 20,000 prime numbers may be used. The first 10,000 of those are considered to be elements belonging to an ordered set $C_{1}'$ and the other 10,000 prime numbers are considered to be elements of a ordered set $C_{2}'$.  Next with $C_{1}'$ and $C_{2}'$, the same operation is carried out as with $C_{1}$ and $C_{2}$. In that way, a new binary string $B'$ is obtained and concatenated to the already existing string. This operation can be repeated with sets of prime numbers  $C_{1}''$ and  $C_{2}''$,  $C_{1}'''$ and $C_{2}'''$, etc.

The procedure to obtain random numbers expressed in base 10 from the last binary string generated is as follows: The first four digits of the string are considered and interpreted as a number in base 2. Hence, for example, the sequences 0000, 0101 and 1001 generate, respectively, digits 0, 5, and 9, in base 10. Each of the sequences 1010, 1011, 1100, 1101, 1110, and 1111 (which generate respectively 10, 11, 12, 13, 14 and 15) is discarded. In this way, the first digit expressed in base 10 is generated, and of course, it may be any of the following: 0, 1, 2,\dots, 9. The second digit in base 10 is generated using the same procedure, from the next sequence of four digits in the binary string; and so on, successively.

\section{Support for the main hypothesis of this paper}

In section 2 mention was made of certain ordered pairs of digits -- $(i, j), i = 0, 1, 2,\dots, 9; j = 0, 1, 2, \dots,  9$ -- such that a) if $i>j$, a 1 is generated as a digit in the binary string from which the random numbers produced by the MRNG are obtained; b) if $i<j$, a 0 is generated as a digit in that string; and c) if $i=j$, then no digit is generated for the string, and the same procedure is followed with the next ordered pair of digits.

The purpose of this procedure is to form a binary string which passes the statistical tests usually applied to determine whether a binary string can be considered random.

The main hypothesis on which this procedure is based is that for anyone who is not familiar with the sequence of all the ordered pairs of digits whose elements are compared, everything takes place as if for each of those comparisons there were the same probability -- $\tfrac{1}{2}$ -- of obtaining a 0 as of obtaining a 1. In other words, if $p(i, j)$ is used to refer to the probability that a pair chosen randomly from the set of random pairs is $(i, j)$, then the following equation must be valid: $p(i, j) = p(j, i)$. If $i \neq j$, the preceding equation implies that, as a result of any of those comparisons made, the probability of obtaining a 1 is equal to that of obtaining a 0.

The best  estimation that can be obtained ``experimentally'' of $p(i,j)$ and $p(j,i)$ are the frequencies $f(i,j)$ and $f(j,i)$, respectively, with which the ordered pairs $(i,j)$ and $(j,i)$ appear. If it is admitted that $p(i, j) = p(j, i)$, then for the corresponding frequencies the following should be true: $f(i,j)\simeq f(j,i)$. This is the verifiable expression of the main hypothesis of this study.

These frequencies were determined, to the fifth decimal place, for 50,000,000 consecutive ordered pairs of the type considered. In table~\ref{t1}, the respective results are given for all of the cases in which $i\neq j$.

\begin{eqnarray*}
f(0,1) = 0.00999 & \simeq & f(1,0) = 0.01001 \\
f(0,2) = 0.01000 & \simeq & f(2,0) = 0.01001 \\
f(0,3) = 0.01001 & \simeq & f(3,0) = 0.00998 \\
f(0,4) = 0.00998 & \simeq & f(4,0) = 0.00999 \\
f(0,5) = 0.01000 & \simeq & f(5,0) = 0.01000 \\
f(0,6) = 0.01000 & \simeq & f(6,0) = 0.00999 \\
f(0,7) = 0.01000 & \simeq & f(7,0) = 0.00996 \\
f(0,8) = 0.00995 & \simeq & f(8,0) = 0.00998 \\
f(0,9) = 0.00998 & \simeq & f(9,0) = 0.00999 \\
f(1,2) = 0.00996 & \simeq & f(2,1) = 0.01001 \\
f(1,3) = 0.00999 & \simeq & f(3,1) = 0.00998 \\
f(1,4) = 0.00998 & \simeq & f(4,1) = 0.01001 \\
f(1,5) = 0.01001 & \simeq & f(5,1) = 0.00998 \\
f(1,6) = 0.01000 & \simeq & f(6,1) = 0.01000 \\
f(1,7) = 0.01002 & \simeq & f(7,1) = 0.00999 \\
\end{eqnarray*}
\begin{eqnarray*}
f(1,8) = 0.00999 & \simeq & f(8,1) = 0.00999 \\
f(1,9) = 0.00999 & \simeq & f(9,1) = 0.00999 \\
f(2,3) = 0.01001 & \simeq & f(3,2) = 0.01000 \\
f(2,4) = 0.01001 & \simeq & f(4,2) = 0.00998 \\
f(2,5) = 0.01000 & \simeq & f(5,2) = 0.00997 \\
f(2,6) = 0.01002 & \simeq & f(6,2) = 0.00999 \\
f(2,7) = 0.01001 & \simeq & f(7,2) = 0.00998 \\
f(2,8) = 0.00998 & \simeq & f(8,2) = 0.00996 \\
f(2,9) = 0.00997 & \simeq & f(9,2) = 0.01003 \\
f(3,4) = 0.01002 & \simeq & f(4,3) = 0.00999 \\
f(3,5) = 0.00999 & \simeq & f(5,3) = 0.01002 \\
f(3,6) = 0.00999 & \simeq & f(6,3) = 0.00999 \\
f(3,7) = 0.01000 & \simeq & f(7,3) = 0.00997 \\
f(3,8) = 0.00997 & \simeq & f(8,3) = 0.01000 \\
f(3,9) = 0.01000 & \simeq & f(9,3) = 0.00998 \\
f(4,5) = 0.00999 & \simeq & f(5,4) = 0.01000 \\
f(4,6) = 0.01002 & \simeq & f(6,4) = 0.01002 \\
f(4,7) = 0.01000 & \simeq & f(7,4) = 0.00996 \\
f(4,8) = 0.01002 & \simeq & f(8,4) = 0.01000 \\
f(4,9) = 0.01001 & \simeq & f(9,4) = 0.01000 \\
f(5,6) = 0.01000 & \simeq & f(6,5) = 0.00999 \\
f(5,7) = 0.00999 & \simeq & f(7,5) = 0.01001 \\
f(5,8) = 0.01001 & \simeq & f(8,5) = 0.01000 \\
f(5,9) = 0.01000 & \simeq & f(9,5) = 0.01000 \\
f(6,7) = 0.00998 & \simeq & f(7,6) = 0.00999 \\
f(6,8) = 0.01000 & \simeq & f(8,6) = 0.00999 \\
f(6,9) = 0.00999 & \simeq & f(9,6) = 0.01000 \\
f(7,8) = 0.01000 & \simeq & f(8,7) = 0.01000 \\
f(7,9) = 0.01002 & \simeq & f(9,7) = 0.00999 \\
f(8,9) = 0.00999 & \simeq & f(9,8) = 0.00999 \end{eqnarray*}\begin{table}[H]\caption{Values of the frequencies $f(i,j)$ -- $i = 0, 1, 2,\dots, 9; j = 0, 1, 2, \dots,  9$ -- for the 90 cases in which $i \neq j$ . Consideration was given to 50,000,000 ordered pairs.}
\label{t1}
\end{table}

In table~\ref{t2}, the frequencies $f(i, j)$ -- $i = 0, 1, 2,\dots, 9; j = 0, 1, 2, \dots,  9$ -- are shown for the 10 cases where $i=j$. It  also can be noted that the frequencies in these cases that do not lead to the generation of digits in the binary string are almost equal.  
\begin{eqnarray*}
f(0,0) = 0.01000 & \qquad & f(5,5) = 0.01000 \\
f(1,1) = 0.00999 & \qquad & f(6,6) = 0.01000 \\
f(2,2) = 0.01001 & \qquad & f(7,7) = 0.00999 \\
f(3,3) = 0.01000 & \qquad & f(8,8) = 0.00999 \\
f(4,4) = 0.01002 & \qquad & f(9,9) = 0.01000 \end{eqnarray*}\begin{table}[H]\caption{Values of the frequencies $f(i,j)$ -- $i = 0, 1, 2,\dots, 9; j = 0, 1, 2, \dots,  9$ -- for the 10 cases in which $i = j$. Consideration was given to 50,000,000 ordered pairs.}
\label{t2}
\end{table}

By observing tables~\ref{t1} and~\ref{t2}, it can be seen that each of the frequencies considered is equal or approximately equal to 0.01. In total, of course, there are 100 frequencies of that type. These tables show that the main hypothesis of this study is confirmed by the values obtained for these frequencies.
\section {Randomness tests and results}

The generator described in section 2 was used to form binary strings which would then be tested for randomness. 

In this section, a discussion will be provided of 1) the statistical tests applied to the binary strings obtained by using the MRNG, and 2) the results of those tests.

Let there be a genuinely random binary string such as that generated by the IRBG. The expression $p(0,0)$ refers to the probability that a dyad (i.e., a sequence of two consecutive bits) selected at random from that string will be $(0, 0)$. Likewise, the expressions $p(0, 1)$, $p(1, 0)$, $p(1, 1)$ will refer to the probabilities that the dyad will be $(0, 1)$, $(1, 0)$ and $(1, 1)$, respectively. For the type of binary string specified the following should be fulfilled: $p(0, 0)$ = $p(0, 1)$ = $p(1, 0)$ = $p(1, 1)$ = $\tfrac{1}{4}$. That is, there is no reason to suppose that in a random binary string a tendency would exist for any one of the four possible dyads to be present, different from the tendency of any of the other three dyads to be present.

The expression $p(0, 0, 0)$ will refer to the probability that a triad (i.e., a sequence of three consecutive bits) which has been randomly chosen from that random binary string, will be (0, 0, 0). Likewise, the expressions $p(0, 0, 1)$, $p(0, 1, 0)$, $p(0, 1, 1)$, $p(1, 0, 0)$, $p(1, 0, 1)$, $p(1, 1, 0)$ and $p(1, 1, 1)$ will refer, respectively, to the probabilities that the triad will be $(0, 0, 1)$, $(0, 1, 0)$, $(0, 1, 1)$, $(1, 0, 0)$, $(1, 0, 1)$, $(1, 1, 0)$ and $(1, 1, 1)$. For the type of binary string specified the following should be fulfilled: $p(0,0,0)$ = $p(0, 0, 1)$ = $p(0, 1, 0)$ = $p(0, 1, 1)$ = $p(1, 0, 0)$ = $p(1, 0, 1)$ = $p(1, 1, 0)$ = $p(1, 1, 1)$ = $\tfrac{1}{8}$. This means that there is no reason to suppose that in a random binary string a tendency would exist for any one of the eight possible triads to be present, different from the tendency of any of the other seven triads to be present.

Clearly, this notion can be generalized for the sixteen possible tetrads, for the thirty-two possible pentads, etc.

Given any two of the consecutive bits in a binary string, there are four possible types of transitions from the first to the second element: from a 0 to a 0, from a 0 to a 1, from a 1 to a 0, and from a 1 to a 1. If the given string is random, the probability that a certain transition (also chosen at random) will be any particular one is equal to the probability that it will be any one of the other three.  That is, there is no reason to suppose that in a random binary string a tendency would exist for any one of the four possible transitions to be present, different from the tendency of any of the other three transitions.

The term ``length of a binary string'' $L_s$ will refer to the number of bits comprising it. Considering the above, for random binary strings of various lengths, one can calculate -- from a theoretical perspective -- the most probable numbers for the different dyads, triads, tetrads, pentads and the transitions that will be present in them. In addition, one can count how many different dyads, triads, tetrads, pentads and transitions are actually present in those binary strings. (The numerical values counted are usually also known as ``values observed''.) It follows quite naturally then that the non-parametric statistical chi-square ($\chi^2$) test may be applied to determine whether, with a given level of significance, the differences found between the expected numerical values based on theoretical considerations and the corresponding numerical values actually observed are significant. In all cases, this test was used with a level of significance of 0.05, with the objective indicated above.

For these cases, the null hypothesis (i.e., the hypothesis that these differences are not significant) is precisely the hypothesis that we wish to prove in this study.

The results of the $\chi^2$ test are shown in table~\ref{t3}, with a level of significance of 0.05 and the degrees of freedom pertinent for the different cases (3 for transitions and dyads, 7 for the triads, 15 for the tetrads and 31 for the pentads) when applied to 1000 binary strings of specified lengths, formed by using the generator being analysed.

\begin{figure}[H]
\centering
\includegraphics[width=4in]{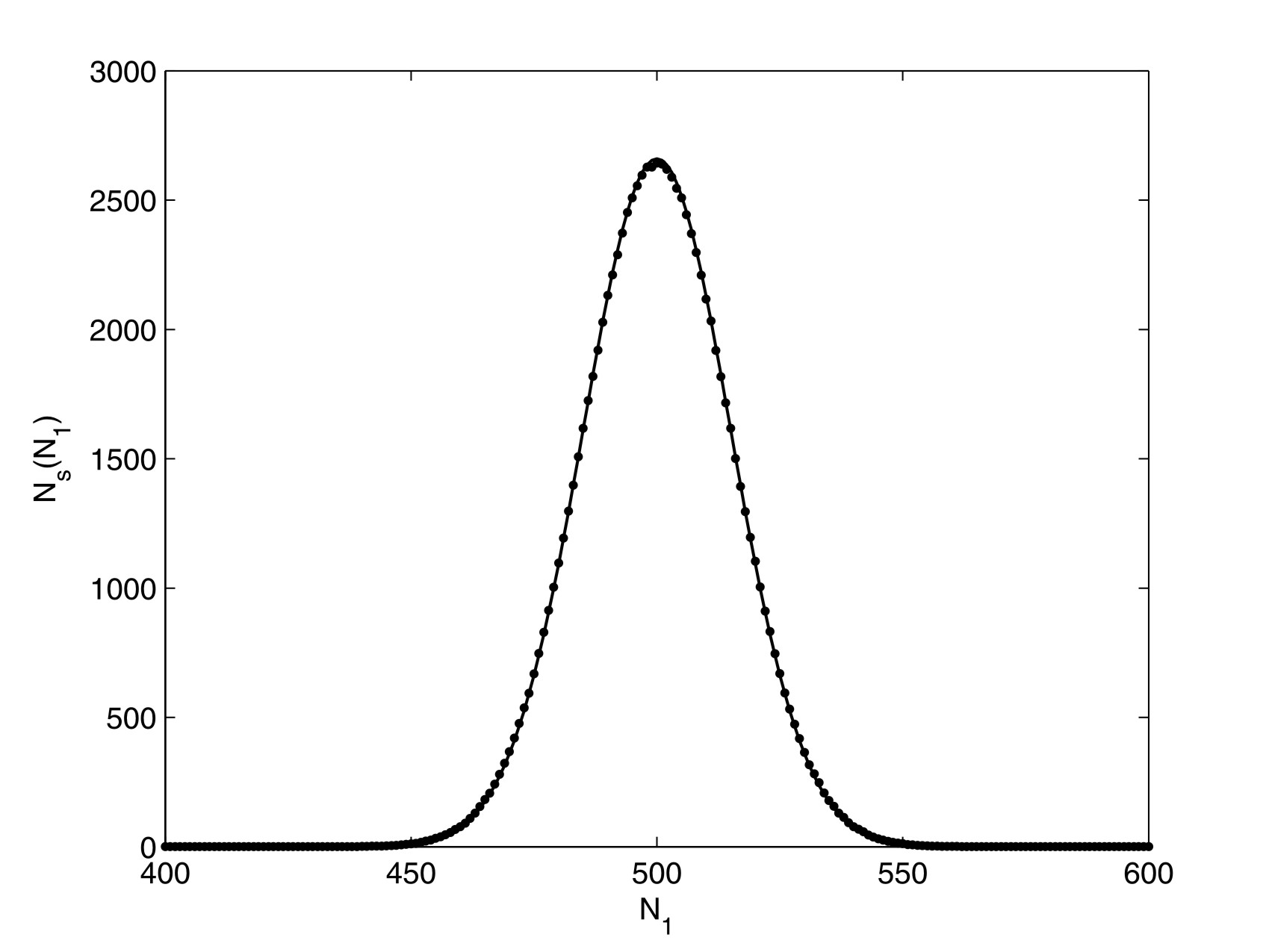}
\caption{Distribution for 100,000 binary strings generated by the MRNG. $N_1$ refers to the number of `ones' generated in each string; $N_s(N_1)$ refers to the number of binary strings as a function of $N_1$.}
\label{f4}
\end{figure}

\begin{table}[H]
\begin{center}
\begin{tabular}{rcccccc}
& \multicolumn{4}{c}{Transitions test} \\  \cline{2-4}
& \multicolumn{4}{c}{$L_s=8001$}\\ \cline{2-4}
&  & MRNG &  \\
\hline 
Passed test & &951&  \\ \hline
Failed test & &49 & \\  \hline
Failed test (\%) & & $4.9\%$& \\  \hline 
\end{tabular}
\end{center}
\end{table}
\begin{table}[H]
\begin{center}
\begin{tabular}{rcccccc}
& \multicolumn{4}{c}{Dyads test} \\  \cline{2-4}
& \multicolumn{4}{c}{$L_s=8000$}\\ \cline{2-4}
&  & MRNG &   \\
\hline 
Passed test & &957 & \\ \hline
Failed test & &43 &  \\  \hline
Failed test (\%) & &$4.3\%$  &  \\ \hline 
\\
& \multicolumn{4}{c}{Triads test} \\  \cline{2-4}
& \multicolumn{4}{c}{$L_s=16000$}\\ \cline{2-4}
&  & MRNG &  \\
\hline 
Passed test & &961 & \\ \hline
Failed test & &39 & \\  \hline
Failed test (\%) &  & $3.9\%$& \\ \hline 
\\
& \multicolumn{4}{c}{Tetrads test} \\  \cline{2-4}
& \multicolumn{4}{c}{$L_s=32000$}\\ \cline{2-4}
&  & MRNG &  \\
\hline 
Passed test & &952 &  \\ \hline
Failed test & &48 & \\  \hline
Failed test (\%)& &$4.8\%$ &  \\ \hline 
\\
& \multicolumn{4}{c}{Pentads test} \\  \cline{2-4}
& \multicolumn{4}{c}{$L_s=64000$}\\ \cline{2-4}
&  & MRNG &   \\
\hline 
Passed test & &952 & \\ \hline
Failed test & &49 &  \\  \hline
Failed test (\%) & & $4.9\% $& \\ \hline 
\end{tabular}
\caption{$\chi^2$ tests}
\label{t3}
\end{center}
\end{table}
The binary strings obtained with the generator considered were subjected  to another test as well. It is based on the following statistical result: The binomial distribution for which each of the two possible events has the same probability of occurring is an excellent approximation to the normal distribution, for very numerous sets of data. Suppose that with that generator, 100,000 binary strings are obtained, where $L_s = 1000$ for each. If a graph is made of the number of binary strings which include exactly a particular number of ``ones'' based on that quantity, a good approximation to a normal curve can be obtained. (Of course, for each of the 100,000 strings the number of ``ones'' included in it could vary from 1 to 1,000.) See this graph in figure~\ref{f4}. (The dotted line represents the binomial distribution and the solid curve the corresponding normal distribution.)

The quality of the approximations (of the binomial distribution considered) to the normal distribution may also be assessed as follows. Let $\mu$ and $\sigma$ be the mean and standard deviation, respectively, of each of these distributions. For both the normal and the binomial distribution obtained by using the $MRNG$, the percentages of the areas under the curves corresponding to those distributions -- for the three specified  intervals -- have been presented in table~\ref{t4}.

\begin{table}[H]
\begin{center}
\begin{tabular}{ccc}
\multicolumn{3}{c}{Percentage of areas under curve} \\ \hline
Interval &    Normal distribution & Binomial distribution considered   \\ \hline
$\mu-\sigma$ to $\mu+\sigma$ &   $68.26\%$  &  $68.27\%$ \\
$\mu-2\sigma$ to $\mu+2\sigma$ & $95.44\%$ &    $95.35\%$ \\
$\mu-3\sigma$ to $\mu+3\sigma$ & $99.74\%$ &    $99.72\%$ \\
\hline
\end{tabular}
\caption{Percentages of area under the curves that correspond to the normal and binomial distributions considered (for the three specified intervals).}
\label{t4}
\end{center}
\end{table}

In a previous article \cite{birnd} we introduced the concept of ``index of randomness'' for $m$-ary strings, where $m=2, 3, 4,\dots$. In particular, this index can be computed for strings such that $m=2$; that is, for binary strings. It is feasible to compute the most probable value of the average of the indexes of randomness of a certain number, such as 25,600, of binary strings that could potentially be generated by the IRBG, such that for each of them $L_s=8$. There are 256 different binary strings with $L_s=8$. For none of them is there a tendency to be generated by the IRBG which is different from the tendency to be generated for any of the other 255 strings. First, then, the randomness index corresponding to each of the aforementioned binary strings is computed. Second, the average is found of the 256 values of the randomness indexes computed. That average is precisely the most likely value sought. In addition, the generator considered here -- the MRNG -- is used to form 25,600 binary strings such that, for each of them, $L_s=8$. With each of the 25,600 binary strings, the following is carried out: First, the randomness index is computed for each of the 25,600 binary strings. Second, the average of the 25,600 randomness indexes already computed is found.

The same approach is used to compute the most probable value of the average of the indexes of randomness of a certain number, such as 6,553,600, of binary strings that could potentially be generated by the IRBG, such that for each of them $L_s=16$. There are 65,536 different binary strings with $L_s=16$. For none of them does there exist a tendency to be generated by the IRBG which is different from the tendency to be generated for any of the other 65,535 strings. First, then, the randomness index corresponding to each of the 65,536 above-mentioned binary strings is computed. Second, the average is found of the 65,536 values of the randomness indexes computed. That average is precisely the most likely value sought. In addition, the generator considered here -- the MRNG -- is used to form 6,553,600 binary strings such that for each of them $L_s=16$. With each of the sets of 6,553,600 binary strings, the following is carried out: First, the randomness index is computed for each of the 6,553,600 binary strings. Second, the average is calculated of the 6,553,600 randomness indexes already computed.

The results are given in Table~\ref{t5}. According to the values of the randomness indexes computed, the generator is of high quality.

\begin{table}[H]
\begin{center}
\begin{tabular}{ccccc}
$L_s=8$ & & & &  \\
& $IRBG$ &  $MRNG$ & \\ \hline

Avg. &0.4238& 0.4237&  \\
Min. &0&0& \\
Max. &0.6466&0.6465&   \\
\hline
\\
$L_s=16$ & &  & &  \\
& $IRBG$  & MRNG &   \\ \hline

Avg. &0.5442 &0.5442&  \\
Min. &0&0&   \\
Max. &0.7113&0.7112&   \\
\hline 
\end{tabular}
\caption{Randomness indexes}
\label{t5}
\end{center}
\end{table}

The preceding results are superior to those that can be obtained by testing (with  the same set of statistical randomness tests)  binary strings generated with two common pseudorandom number generators, such as JAVA and C++.  These results are also up to par, quality-wise, with those of the statistical tests performed on the binary string obtained by another random number generator previously developed by the authors \cite{bhrng}.

\section{Discussion and perspectives}

In S. B. Volchan's article ``What is a random sequence?''\cite{bvolchan}, it is possible to find the answers which authors with different theoretical conceptions have given to that very question. (p. 48) At present the best known notion of randomness is perhaps that of Solomonoff, Kolmogorov and Chaitin. This approach -- that of randomness as incompressiblity -- applied to binary strings, establishes that a binary string is random if there is no computer program which can generate that string and whose length, in bits, is less than the length of that string (also characterized by the number of digits, each of which can be a 0 or a 1, comprising it). This approach is quite developed, and has resulted in Gregory Chaitin's Algorithmic Information Theory (AIT) (\cite{bc1},\cite{bc2},\cite{bc3}). The AIT has been very useful in the field of metamathematics, but does not make it possible to generate random sequences of digits, even though the concept of ``random'' is considered in the sense accepted by that theory.

The binary string, made up of digits, whose generation technique was specified above, is not random according to the AIT, because a computer program does exist (that developed by the authors of this paper), which has generated it and whose length (in bits) is less than the length (also in bits) of this string. It $is$ random, however,  according to the criterion usually found in statistics to evaluate whether a binary string is random: that of passing a certain ``battery'' of randomness tests, currently accepted as suitable.

In fields such as statistics, computer simulation, and cryptography, it may be useful to consider approaches which differ from the dichotomy ``random string--regular (or non-random) string''. In this regard, it is pertinent to quote from \cite{birnd}:

``S. Pincus and B. H. Singer have developed an approach which varies, partially, from the dichotomic conception `random strings--nonrandom strings' \cite{bpincus}. According to these authors, `...there appears to be a critical need to develop a (computable) formulation of `closer to random', to grade the large class of decidedly nonrandom sequences'. They accept the existence of a `large class of decidedly nonrandom sequences' and consider it appropriate to distinguish in that class between various levels of nearness to truly random sequences. This approach is different from ours, where for binary strings of a given length, there are only two of them (and not one `large class') which can be considered definitely regular (or nonrandom): that comprised only of ones and that comprised only of zeros. For these two strings, the index of regularity is equal to 1 and the index of randomness is equal to 0. The remaining strings are not considered to be either definitely random or definitely regular. For each of these strings, it is feasible to compute an index of regularity ($i_{reg}$) and an index of randomness ($i_{rnd}$), such that $i_{reg} + i_{rnd} = 1$. It is obvious that, following our approach, the higher the regularity index of a binary string, the lower its randomness index, and vice versa. With this approach there are as many classes of degrees of randomness as there are different indexes of randomness computed. The higher the index of randomness of a string, the more it can be considered random. Of course, for $m$-ary strings (with $m = 3, 4, 5,\dots$) of a given length, there are $m$ strings for each of which $i_{reg} = 1$ and $i_{rnd} = 0$.'' (p. 57--58)

Information will be provided here concerning what was meant by ``in detail'' in section 1, with regard to how much is known about the procedure used to generate  the binary string from which the random numbers are obtained. That expression means that precisely this information is known:
\\1) with how many digits the decimal parts of the roots of the prime numbers are obtained;
\\2) how many digits (between 40 and 100) are discarded upon using, as specified above, the decimal parts of the different roots of the prime numbers; and 
\\3) what permutations $P_{1}$, $P_{2}$, $P_{3}$,$\dots$,$P_{10,000}$ were applied to the prime numbers belonging to set $C_{2}$.

If one wants to use the MRNG not only in the field of statistics -- for example, to select random samples from a given population -- and in the field of computer simulation, but also in that of cryptography, that information must not be provided to those who are  not authorized to decipher messages encrypted by using the binary string obtained. In effect, having not only the knowledge about the procedure described in section 2 but also that information makes it possible to reobtain that binary string with no ambiguity whatsoever, and therefore, the messages which were encrypted by using it.

The Mathematical Random Number Generator (MRNG) facility may be found on the webpage of the ``Applied Mathematics \& Computer Simulation Group'' -- www.appliedmathgroup.org -- and is available for use.  Thus, anyone interested in verifying the quality of this generator may test the binary string generated by it as well as the random numbers obtained from it, using suitable statistical tests. Thus for example, in the second section, ``2. Random Number Generation Tests'', of a document published by the National Institute of Standards and Technology (NIST) \cite{bnist}, one can find the ``NIST Test Suite'' \textbf{--} a statistical package consisting of 15 tests developed to test the randomness of (arbitrarily long) binary sequences produced by either hardware- or software-based cryptographic random or pseudorandom number generators.

It should be emphasized that if an infinite sequence of digits were produced by a generator of the type described here, it would $not$ be a periodic sequence, like an infinite sequence would be if it were produced by a congruential generator, such as that used in Java and C++.

The MRNG is easy to implement, and it has many acceptable variants. Selected variants which already have been designed will be described elsewhere.

\end{document}